\documentclass[12pt]{article}  
\usepackage{amsmath}
\usepackage{hyperref}
\usepackage{amsfonts}
\usepackage{amssymb}
\usepackage{graphics}
\usepackage{color}
\usepackage[pdftex]{graphicx}

\def\bee{\begin{equation}}
\def\eee{\end{equation}}

\def\OO{\mathcal{O}}
\def\hksqrt{\mathpalette\DHLhksqrt}
\def\DHLhksqrt#1#2{\setbox0=\hbox{$#1\sqrt{#2\,}$}\dimen0=\ht0
\advance\dimen0-0.2\ht0
\setbox2=\hbox{\vrule height\ht0 depth -\dimen0}%
{\box0\lower0.4pt\box2}}
\begin{document}

\thispagestyle{empty}
\bigskip
\centerline{    }
\centerline{\Large\bf A Note on the Andrica Conjecture}
\bigskip\bigskip
\centerline{\large\sl Marek Wolf}

\begin{center}
e-mail: primes7@o2.pl
\end{center}

\bigskip\bigskip

\begin{center}
{\bf Abstract}\\
\end{center}

\begin{minipage}{12.8cm}
We derive heuristically the approximate formula for the difference
$\sqrt{p_{n+1}} - \sqrt{p_n}$, where $p_n$ is the n-th prime. We find perfect
agreement between this formula  and the available data  from
the list of maximal  gaps between consecutive primes.
\end{minipage}

\bigskip\bigskip

\bigskip\bigskip

\bibliographystyle{abbrv}

\section{Introduction}

The Andrica conjecture \cite{Andrica} (see also \cite[p.21]{Guy} and
\cite[p. 191]{Ribenboim}) states that the inequality:
\bee
A_n \equiv \sqrt{p_{n+1}} - \sqrt{p_n} < 1
\label{Andrica-ineq}
\eee
where $p_n$ is the $n$-th prime number, holds for all $n$.  Despite its simplicity
it remains unproved. In the Table I we give a few first values of $A_n$ and
in Table II the values of $A_n$ are sorted in descending order.

We have
\bee
\sqrt{p_{n+1}} - \sqrt{p_n}=\frac{p_{n+1} - p_n}{\sqrt{p_{n+1}} + \sqrt{p_n}}<
\frac{d_n}{2\sqrt{p_n}}
\label{nierownosc}
\eee
From this we see that the growth rate of the form $d_n = \OO(p_n^\theta)$ with
$\theta<1/2$ will suffice for the proof of (\ref{Andrica-ineq}).
Unfortunately all  values of $\theta$ proved in the past are larger than $1/2$.
A few results with  $\theta$ closest to 1/2 are: M. Huxley: $\theta>7/12$ \cite{Huxley},
the  result of C.J. Mozzochi  \cite{Mozzochi1986} $\theta = {1051\over 1920}$,
S. Lou and Q. Yao  obtained $\theta=6/11$ \cite{Lou1992} and recently R.C. Baker
G. Harman and J. Pintz  \cite{Baker-et-al} have improved it to
$\theta=21/40$ what remains currently the best unconditional result.
For a review of results on $\theta$ see \cite{Pintz-Landau}.
The best estimation  for $d_n$ obtained by H. Cramer under the assumption of
the Riemann Hypothesis \cite{Cramer1920}
\bee
d_n=\mathcal{O}(\sqrt{p_n}\log(p_n))
\eee
also does not suffice to prove the Andrica conjecture.

\vskip 0.4cm
\begin{center}
{\sf TABLE {\bf I}}\\
\bigskip
\begin{tabular}{||c|c|c|c||c|c|c|c||} \hline
$p_n$ & $ p_{n+1} $ & $ d_n $ & $\sqrt{p_{n+1}}-\sqrt{p_n}$ & $p_n$ & $ p_{n+1} $ & $ d_n $ & $\sqrt{p_{n+1}}-\sqrt{p_n}$ \\ \hline
     2  &       3 &      1 &         0.317837245  &   41  &      43 &      2 &         0.154314287    \\ \hline
     3  &       5 &      2 &         0.504017170  &   43  &      47 &      4 &         0.298216076    \\ \hline
     5  &       7 &      2 &         0.409683334  &   47  &      53 &      6 &         0.424455289    \\ \hline
     7  &      11 &      4 &         0.670873479  &   53  &      59 &      6 &         0.401035859    \\ \hline
    11  &      13 &      2 &         0.288926485  &   59  &      61 &      2 &         0.129103928    \\ \hline
    13  &      17 &      4 &         0.517554350  &    61  &      67 &      6 &         0.375103096   \\ \hline
    17  &      19 &      2 &         0.235793318  &    67  &      71 &      4 &         0.240797001   \\ \hline
    19  &      23 &      4 &         0.436932580  &    71  &      73 &      2 &         0.117853972   \\ \hline
    23  &      29 &      6 &         0.589333284  &    73  &      79 &      6 &         0.344190672   \\ \hline
    29  &      31 &      2 &         0.182599556  &    79  &      83 &      4 &         0.222239162   \\ \hline
    31  &      37 &      6 &         0.514998167  &    83  &      89 &      6 &         0.323547553   \\ \hline
    37  &      41 &      4 &         0.320361707  &    89  &      97 &      8 &         0.414876670   \\ \hline
    41  &      43 &      2 &         0.154314287  &    97  &     101 &      4 &         0.201017819   \\ \hline
    43  &      47 &      4 &         0.298216076  &   101  &     103 &      2 &         0.099015944   \\ \hline
    47  &      53 &      6 &         0.424455289  &   103  &     107 &      4 &         0.195188868   \\ \hline
    53  &      59 &      6 &         0.401035859  &   107  &     109 &      2 &         0.096226076   \\ \hline
    59  &      61 &      2 &         0.129103928  &   109  &     113 &      4 &         0.189839304   \\ \hline
\end{tabular} \\
\end{center}
\vskip 0.4cm

For twins primes $p_{n+1}=p_n+2$ there is no problem
with (\ref{Andrica-ineq}) and in general for short gaps $d_n=p_{n+1} - p_n$
between consecutive primes the inequality (\ref{Andrica-ineq}) will be satisfied.
The Andrica conjecture can be violated only by  extremely large gaps between
consecutive primes.  Let $G(x)$ denote the largest gap between
consecutive primes smaller than $x$:
\bee
G(x)={\max_{p_n, p_{n-1}<x}} (p_n-p_{n-1}).
\eee
Let us denote the pair of primes $<x$  comprising the largest gap $G(x)$
by $p_{L+1}(x)$ and $p_L(x)$, hence we have
\bee
G(x)=p_{L+1}(x)- p_L(x).
\eee
Thus we will concentrate on the values of the difference appearing in
(\ref{Andrica-ineq}) corresponding to the largest gaps and let us
introduce the function:
\bee
R(x)=\sqrt{p_{L+1}(x)}- \sqrt{p_L(x)}
\eee
Then we have:
\bee
A_n\leq R(p_n).
\eee
The largest values of $A_n$ will be reached at the largest gaps $G(x)$ between consecutive
primes  below a given bound $x$. In  \cite{Borcherds-Krakow}, \cite{Wolf-conj}
we have given the heuristic arguments that  $G(x)$ can be expressed 
directly by $\pi(x)$ --- the number of primes $<x$:
\bee
G(x) \sim {x \over \pi(x)} (2\log \pi(x) -\log(x) +c'),
\label{G-max}
\eee
where $c'$ is expressed by the twin constant $C_2$:
\bee
c'=\log(C_2)=0.27787688\ldots , ~~~~C_2 \equiv  2\prod_{p > 2} \biggl( 1 - {1 \over (p - 1)^2}\biggr) =
1.32032363169\ldots
\eee
For the Gauss approximation $\pi(x)\sim x/\log(x)$  the following dependence follows:
\bee
G(x) \sim \log(x) (\log(x)-2\log\log(x) +\log(c'))
\label{d_max2}
\eee
and for large $x$ it passes into the Cramer \cite{Cramer} conjecture:
\bee
G(x)\sim \log^2(x).
\label{eCramer}
\eee

\vskip 0.4cm
\begin{center}
{\sf TABLE {\bf II}}\\
\bigskip
\begin{tabular}{||c|c|c|c|c||} \hline
$ n $ & $p_n$ & $ p_{n+1} $ & $ d_n $ & $\sqrt{p_{n+1}}-\sqrt{p_n} $\\ \hline
4  &       7 & 11  & 4 & 0.6708735 \\ \hline
30  &      113 & 127  & 14 & 0.6392819 \\ \hline
9  &       23 & 29  & 6 & 0.5893333 \\ \hline
6  &       13 & 17  & 4 & 0.5175544 \\ \hline
11  &      31 & 37  & 6 & 0.5149982 \\ \hline
2  &       3 & 5  & 2 & 0.5040172 \\ \hline
8  &       19 & 23  & 4 & 0.4369326 \\ \hline
15  &      47 & 53  & 6 & 0.4244553 \\ \hline
46  &      199 & 211  & 12 & 0.4191031 \\ \hline
34  &      139 & 149  & 10 & 0.4167295 \\ \hline
\end{tabular} \\
\end{center}
\vskip 0.4cm

A. Granville argued  \cite{Granville} that the actual $G(x)$ can be larger
than that given by (\ref{eCramer}), namely he claims  that
there are infinitely many pairs of primes $p_n, p_{n+1}$ for which:
\bee
p_{n+1} - p_n = G(p_n)> 2 e^{-\gamma}\log^2(p_n) = 1.12292\ldots \log^2(p_n).
\label{Granville-Cramer}
\eee

For a given gap $d$ the largest value of the difference $\sqrt{p+d}-\sqrt{p}$
will appear at the first appearance of this gap:
each next pair $(p', p'+d)$ of consecutive primes separated by $d$ will
produce smaller difference (see (\ref{nierownosc})):
\bee
\sqrt{p'+d}-\sqrt{p'}<\sqrt{p+d}-\sqrt{p}.
\eee
Hence we have to focus our attention on the first occurrences of gaps.
In  \cite{Wolf-first-d} we have given heuristic
arguments that the gap $d$ should appear for the first time after the prime
$p_f(d)$ given by
\bee
p_f(d) \sim \sqrt{d} e^{\sqrt{d}}.
\eee
We calculate

\bee
\begin{split}
\sqrt{p_f(d)+d}-\sqrt{p_f(d)}=\sqrt{\sqrt{d} e^{\sqrt{d}}+d}-\sqrt{\sqrt{d} e^{\sqrt{d}}}=\\
\sqrt{\sqrt{d} e^{\sqrt{d}}}\Big(\hksqrt{1+\frac{d}{\sqrt{d} e^{\sqrt{d}}}}-1\Big)=
\frac{1}{2}d^{\frac{\small 3}{\tiny 4}}e^{-\frac{1}{2}\sqrt{d}}+\ldots
\end{split}
\label{wykladki}
\eee
Substituting here for $d$ the maximal gap $G(x)$ given by (\ref{G-max})
we obtain the approximate formula for $R(x)$:
\bee
R(x)=\frac{1}{2}G(x)^{3/4}e^{-\frac{1}{2}\sqrt{G(x)}}+error~~term.
\label{main-formula}
\eee


The comparison with real data is given in Figure 1. The lists of known maximal
gaps between consecutive primes can be found at http://www.trnicely.net
and  http://www.ieeta.pt/$\sim$tos/gaps.html. The largest known  gap 1476
between consecutive primes follows the prime  $1425172824437699411=
1.42\ldots\times10^{18}$.

The maximum of the function $\frac{1}{2}x^{\frac{\small 3}{\tiny 4}}e^{-\frac{1}{2}\sqrt{x}}$
is reached at $x=9$ and has the value $0.579709161122$. The maximal value of $A_n$ is 
$0.6708735\ldots$ for $d=4$ and second value is $0.6392819\ldots$ for $d=14$. 
Let us remark that $d=9$ is exactly in the middle between $4$ and 14.

\begin{figure}
\vspace{-3.1cm}
\hspace{-2.5cm}
\begin{center}
\includegraphics[height=0.7\textheight,angle=90]{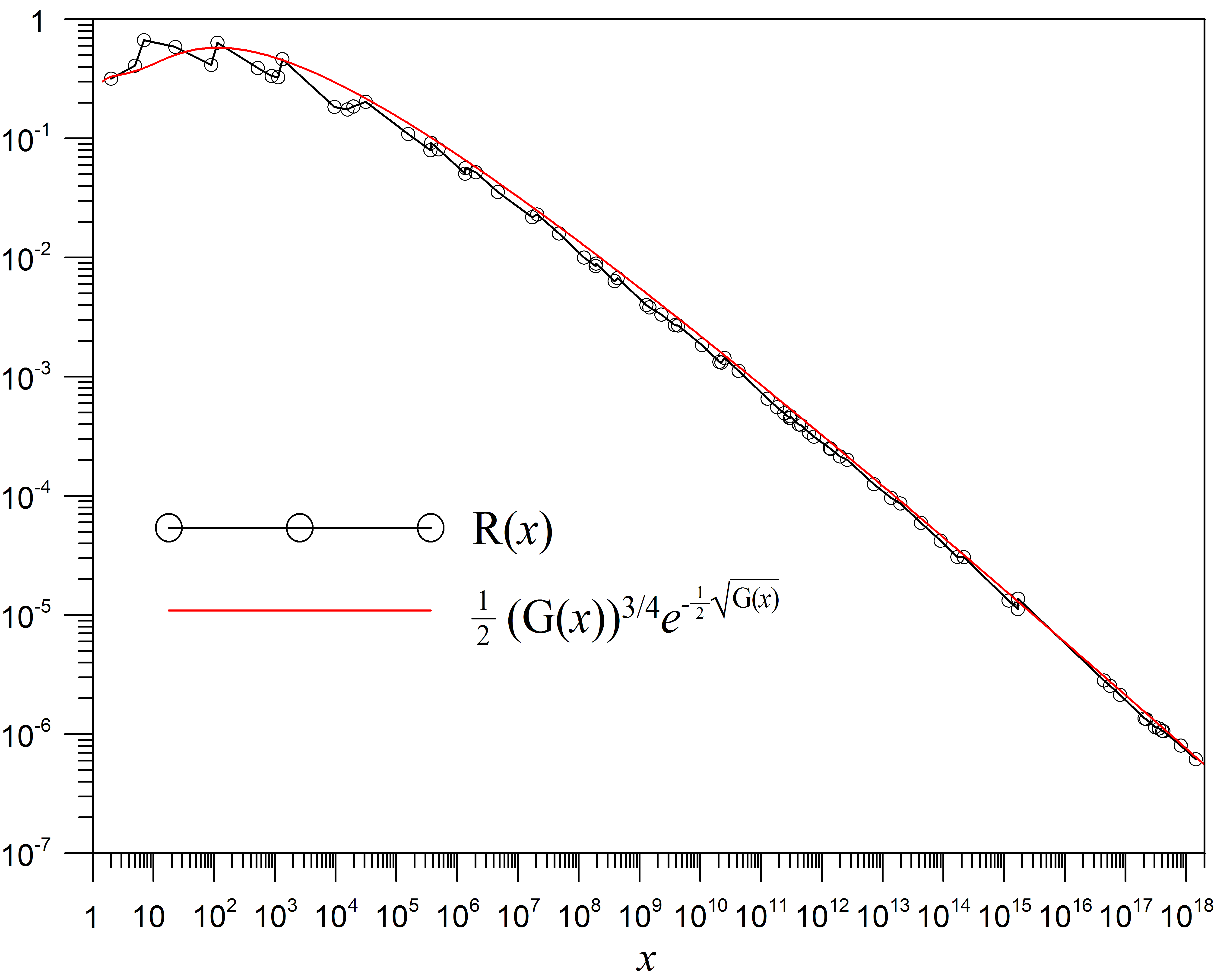} \\
\vspace{0.0cm} {\small Fig.1 The plot of $R(x)$ and approximation to it
given by ({main-formula}). The are 75 maximal gaps available currently and
hence there are 75 circles in the plot of $R(x)$.}
\end{center}
\end{figure}

Because in (\ref{main-formula}) $R(x)$ contains exponential of $\sqrt{G(x)}$ it is
very  sensitive to the form of $G(x)$.  The substitution  $G(x)=\log^2(x)$
leads to the form:
\bee
R(x)=\frac{\log^{3/2}(x)}{2\sqrt{x}}.
\eee
This form of $R(x)$ is plotted in Fig.2 in red. In \cite{Shanks1964} D. Shanks has
given for $p_f(d)$ the expression
\bee
p_f(d) \sim e^{\sqrt{d}}.
\eee
This leads to the expression
\bee
\sqrt{p_f(d)+d}-\sqrt{p_f(d)}=\frac{1}{2} d e^{-\frac{1}{2}\sqrt{d}}
\eee
instead of (\ref{wykladki}). Substitution here for $d$ the form (\ref{d_max2})
leads to the curve plotted in Fig.2 in green.

\begin{figure}
\vspace{-1.1cm}
\hspace{-1.5cm}
\begin{center}
\includegraphics[height=0.4\textheight,angle=0]{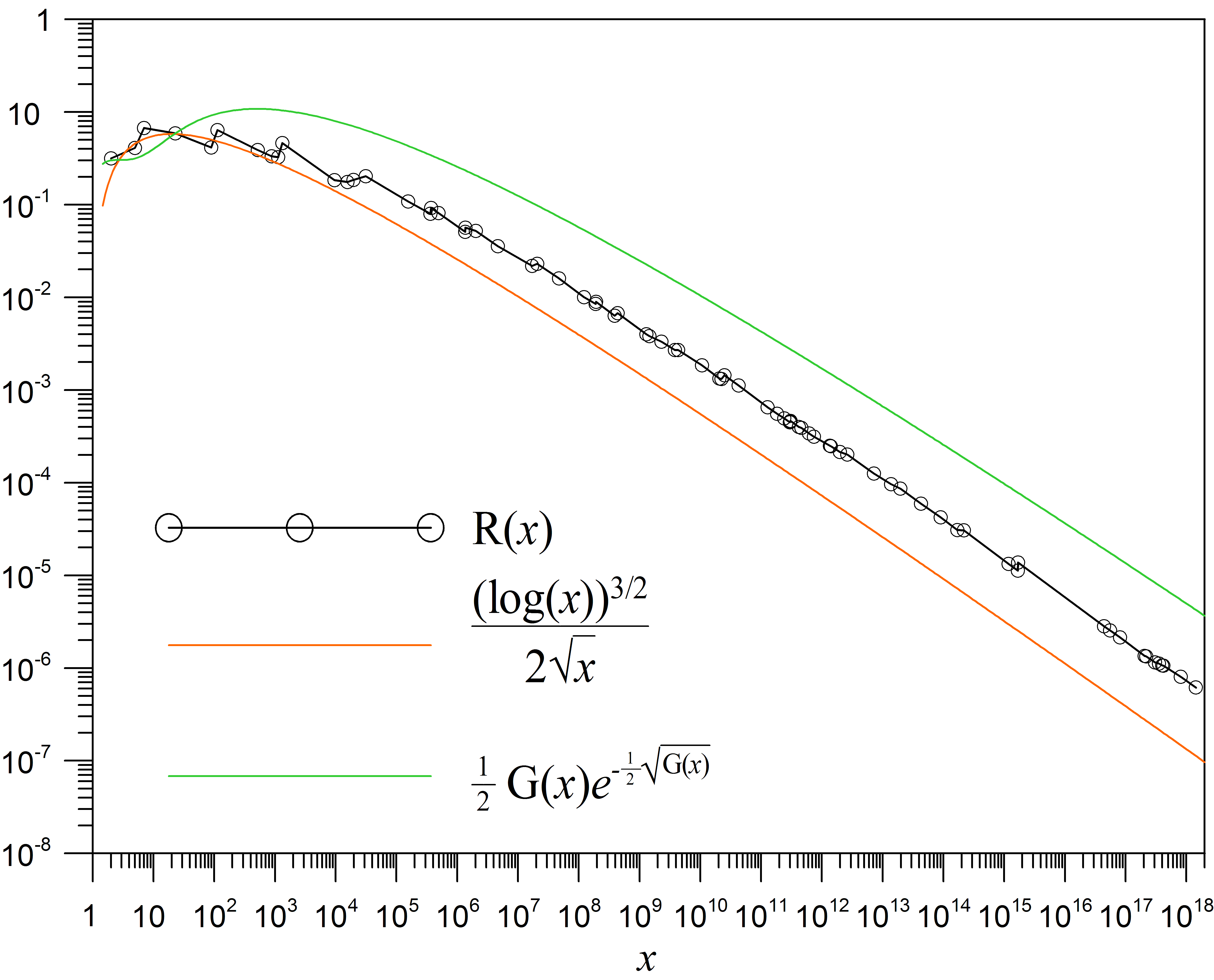} \\
\vspace{0.0cm} {\small Fig.1 The plot of $R(x)$ and approximation to it
given by $\frac{\log^{3/2}(x)}{2\sqrt{x}}$ (red) and approximation of $R(x)$
obtained from the Shanks conjecture for $p_f(d)$ (green). }
\end{center}
\end{figure}

Finally let us remark, that  from the above analysis it follows, that
\bee
\lim_{n \to \infty} \sqrt{p_{n+1}} - \sqrt{p_n} =0
\eee
The above  limit  was mentioned on p. 61 in \cite{Golomb1976} as a difficult problem
(yet unsolved).


\begin{thebibliography}{10}

\bibitem{Andrica}
D.~Andrica.
\newblock Note on a conjecture in prime number theory.
\newblock {\em Studia Univ. Babes-Bolyai Math.}, 31:44--48, 1986.

\bibitem{Baker-et-al}
R.~C. Baker, G.~Harman, and P.~J.
\newblock The difference between consecutive primes, {II}.
\newblock {\em Proc. London Math. Soc.}, 83:532–562, 2001.

\bibitem{Cramer1920}
H.~Cramer.
\newblock Some theorems concerning prime numbers.
\newblock {\em Arkiv f. Math. Astr. Fys.}, 15:1–--33, 1920.

\bibitem{Cramer}
H.~Cramer.
\newblock On the order of magnitude of difference between consecutive prime
  numbers.
\newblock {\em Acta Arith.}, II:23--46, 1937.

\bibitem{Golomb1976}
S.~W. Golomb.
\newblock Problem {E}2506: Limits of differences of square roots.
\newblock {\em Amer. Math. Monthly}, 83:60--61, 1976.

\bibitem{Granville}
A.~Granville.
\newblock Harald {C}ramer and the distribution of prime numbers.
\newblock {\em Scandanavian Actuarial J.}, 1:12--28, 1995.

\bibitem{Guy}
R.~K. Guy.
\newblock {\em Unsolved Problems in Number Theory}.
\newblock Springer-Verlag, 2nd ed. New York, 1994.

\bibitem{Huxley}
M.~Huxley.
\newblock An application of the {F}ouvry-{I}waniec theorem.
\newblock {\em Acta Arithemtica}, 43:441--443, 1984.

\bibitem{Lou1992}
S.~Lou and Q.~Yao.
\newblock A {C}hebyshev's type of prime number theorem in a short interval. ii.
\newblock {\em Hardy-Ramanujan J.}, 15:1--33, 1992.

\bibitem{Mozzochi1986}
C.~Mozzochi.
\newblock On the difference between consecutive primes.
\newblock {\em Journal Number Theory}, 24:181--187, 1986.

\bibitem{Pintz-Landau}
J.~Pintz.
\newblock Landau's problems on primes.
\newblock {\em Journal de théorie des nombres de Bordeaux}, 21:357--404, 2009.

\bibitem{Ribenboim}
P.~Ribenboim.
\newblock {\em {The Little Book of Big Primes}}.
\newblock 2ed., Springer, 2004.

\bibitem{Shanks1964}
D.~Shanks.
\newblock On maximal gaps between successive primes.
\newblock {\em Mathematics of Computation}, 18:646--651, 1964.

\bibitem{Wolf-conj}
M.~Wolf.
\newblock Some conjecturees on the gaps between consecutive primes, 1995.
\newblock preprint IFTUWr 894//95, September 1995, available from
  http://www.ift.uni.wroc.pl/$\sim$mwolf/conjectures.ps.gz.

\bibitem{Borcherds-Krakow}
M.~Wolf.
\newblock Unexpected regularities in the distribution of prime numbers.
\newblock In P.~et~al, editor, {\em 8th Joint EPS-APS Int.Conf. Physics
  Computing'96, Kraków, 1996}, pages 361--367, 1996.

\bibitem{Wolf-first-d}
M.~Wolf.
\newblock First occurence of a given gap between consecutive primes, 1997.
\newblock preprint IFTUWr 911//97, April 1997, available from
  http://www.ift.uni.wroc.pl/$\sim$mwolf/firstocc.pdf.

\end{thebibliography}
\end{document}